\documentstyle[12pt]{article}
\newtheorem{theorem}{Theorem}[section]
\newtheorem{definition}[theorem]{Definition}
\newtheorem{lemma}[theorem]{Lemma}
\newtheorem{corollary}[theorem]{Corollary}
\newtheorem{proposition}[theorem]{Proposition}
\begin{document}

\centerline{\large \bf Uniform boundedness for rational points}

\smallskip

\leftline{Patricia L. Pacelli}
\leftline{Department of Mathematics, Boston University}
\leftline{111 Cummington Street, Boston, MA 02215}
\leftline{{\tt pacelli@math.bu.edu}}

\smallskip

\noindent \today

\noindent This is a preliminary version

\section{Introduction}\label{intro}

According to a famous conjecture of S. Lang's, if $K$ is a number
field, then the set of $K$-rational points of any variety of general
type defined over $K$ is not dense in the Zariski topology.  In a
recent paper, entitled {\it Uniformity of rational points}
(\cite{chm}), L. Caporaso, J. Harris, and B. Mazur show that this
conjecture implies the existence of a uniform bound on the number of
$K$-rational points over all smooth curves of genus $g$ defined over
$K$, for some fixed $g\geq 2$.  Their bound depends on the genus $g$ and
on the number field $K$.

D. Abramovich has proved an extension of this result.  In his paper,
\cite{a}, he proves that,
assuming Lang's conjecture, the bound $B(K,g)$ of \cite{chm} remains
bounded as $K$ varies over all quadratic number fields, or as $K$
varies over all quadratic extensions of a fixed number field.

It is this result of Abramovich's that we shall generalize in this
paper.  We will prove that given a number field $K$, Lang's conjecture
implies the existence of a uniform bound on the number of $L$-rational
points over all smooth curves of a fixed genus $g>1$ defined over $L$,
as $L$ varies over all extensions of $K$ of degree $d$ for any
positive integer $d$.  This bound will depend on $K$, $d$, and $g$,
but is independent of the actual number field $L$.

\begin{theorem}\label{mainth} Assume that Lang's conjecture regarding
varieties of general type is true. Let $g\geq 2$ and $d\geq 1$ be
integers, and let $K$ be a number field.  Then there exists an integer
$B_K(d,g)$, which, for a given $K$ depends only $d$ and $g$, such that for any
extension $L$ of $K$ of degree $d$, and any curve $C$ of genus $g$
defined over $L$, it follows that $$\#C(L)\leq B_K(d,g).$$
\end{theorem}

By letting $K={\bf Q}$ we have the following:

\begin{corollary}\label{cor} Assume Lang's
conjecture is true. Let $g\geq 2$ and
$d\geq 1$ be integers.
Then there exists a bound $B(d,g)$, depending only on $d$ and $g$,
such that for any number field $L$ of degree $d$, and for any curve $C$
of genus $g$ defined over $L$, it follows
that
$$\#C(L)\leq B(d,g).$$
\end{corollary}

Since any extension $L$ of $K$ of degree $d$ is a number field of
some fixed degree $d'$, Theorem \ref{mainth} and Corollary \ref{cor}
are equivalent.  We
state Theorem \ref{mainth} separately, however, as it might be interesting
to study the dependence of the bound on $K$ in later work, perhaps
to see what happens assuming Lang's so-called strong conjecture.

\subsection{Acknowledgements} I am extremely grateful
to Dan Abramovich for introducing me to this interesting problem.  His
guidance and assistance have been invaluable.  I would also like to
thank Emma Previato and Lucia Caporaso for listening to presentations
of this material, and Joseph Silverman for giving me advice and
comments.  Finally, I owe a tremendous amount of thanks to Glenn
Stevens, for both the mathematics he has taught and the encouragement
he has provided.

\section{Definitions, Notation, and Ideas}

Lang's conjecture concerns varieties of general type.  Recall the
definitions:

\begin{definition}\label{def1} A line bundle $L$ on a variety $Y$ is said to be
{\em big} if for high values of $k$, $H^0(Y,L)$ has enough sections to induce
a birational map to projective space.
\end{definition}

\begin{definition} A smooth projective variety $Y$ is of {\em general type} if
its dualizing sheaf $\omega_Y$ is big.  An arbitrary
projective variety is of general type if a desingularization of it is.
\end{definition}

The main idea we will be working with is that given a family of curves
$X\to B$, we can study the symmetric $d$-th power of the
$n$-th fibered product (for some sufficiently large $n$),
$$Sym^d(X^n_B)={(X_B^n)}^d/S_d.$$

For our purposes we will take the family $X\to B$ that we
work with to be a {\it tautological family}.
This is a family
$X\to B$ of stable curves along with a finite surjective map
$\phi\colon B\to \overline{M_g}$, where $\overline{M_g}$ is
the
moduli space of smooth stable curves of genus $g$.  The map $\phi$
is assumed to have the property that
$\phi(b) = [X_b]$ for all $b\in B$.  See \cite{chm}, \S 5.1 for a proof that
such a family exists.  The reason we work with a tautological
family is quite simple; our theorem asserts that there is a bound on the
number of rational points on any smooth curve of genus $g$ - hence we want a
family in which every such curve appears as a fiber, possibly after a
field extension of bounded degree.

Assume $d\geq1$ and $g\geq2$ to be fixed integers throughout
the remainder of the paper, where $d$ represents the degree of an extension
of number fields $K\subseteq L$, and $g$ represents the genus of the curves we
will be looking at.  Following the notation in \cite{a}, write
$$Y_n = Sym^d(X_B^n).$$

Let $L$ be an extension of $K$ of degree $d$, and let $\sigma_1$,
\ldots,$\sigma_n$ be the $d$ embeddings of $L$ in $\overline{K}$
fixing $K$.  Then
if $b\in B(L)$, and $(P_1, \ldots, P_n) \in X_b(L)$, we obtain a
$K$-rational point $y_{(P_1, \ldots, P_n)}$ of $Y_n$:
$$
y_{(P_1,\ldots, P_n)} = \{(P_1, \ldots, P_n)^{\sigma_1}, \ldots,
 (P_1, \ldots, P_n)^{\sigma_d}\}
$$

\begin{definition} For a variety $V$ defined over $K$ and a field extension
$K\subset L$, let $V(L,K)$ be the set of those points lying
on the variety $V$ which are defined over $L$ but are not
defined over $K$ nor any other intermediate field between
$K$ and $L$.
\end{definition}

\begin{definition} Let $m\geq n$.  Again following \cite{a}, we call a point
$y_{(P_1,\ldots,P_n)}\in Y_n(K)$
{\it m-prolongable}
if there exists a number field $L$, $[L\colon K] = d$, and a
{\it prolongation} $y_{(P_1,\ldots,P_m)}\in Y_m(K)$ such that
$P_i\in X(L,K)$ for all $1\leq i\leq m$, and if
$1\leq i\not= j\leq m$ then $P_i\not= P_j.$  In other words,
$y_{(P_1,\ldots,P_n)}\in Y_n$ is m-prolongable if the n points
$P_1,\ldots,P_n$ can be extended to m distinct
points $P_1,\ldots P_m \in X(L,K)$, producing a point
$y_{(P_1,\ldots,P_m)}\in Y_m.$
\end{definition}

We will call $E_n^{(m)}$ the set of $m$-prolongable points in $Y_n$, and
will denote by $F_n^{(m)}$ the Zariski closure $\overline{E_n^{(m)}}$.

\begin{lemma} For all $n\geq 1$, there exists an integer
$m(n)$ such that for all positive integers $k$,
$$F_n^{(m(n)+k)}=F_n^{(m(n))}.$$
\end{lemma}

{\bf Proof:} Since $F_n^{(m+1)}\subseteq  F_n^{(m)}$, the $F_n^{(m)}$'s
form a decreasing sequence of closed sets in the noetherian space $Y_n$,
which must eventually stop.  So for all $n$ there is an integer $m(n)$
such that $F_n^{(m(n)+k)}=F_n^{(m(n))}$ for all positive integers $k$.
Q.E.D.

To ease notation we will write $F_n$ for $F_n^{(m(n))}$.

Our goal is to show that each $F_n$ is empty.  For each $n$, $F_n$ is
the closure of those $n$-tuples of distinct points defined over $L$
which can be extended to arbitrarily long $m$-tuples of distinct
points over the same field.  If we show that each $F_n$ is empty, this
shows that there must be a bound of the desired type.

We want to look at some of the properties of the $F_n$'s, and study the
maps between them.  The ideas in the remainder of this section are
generalizations of those
found in \cite{a}.

\begin{lemma} Suppose $n'>n$ and $I\subseteq \{1, 2, \ldots, n'\}$ is an
$n$-tuple.  Then the projection $\pi_I :F_{n'}\to F_n$ is surjective.
\end{lemma}

{\bf Proof:} This is clearly true for the $E_n^{(m)}$'s by definition.
Thus the same holds for the $F_n$'s.  Q.E.D.

We have a natural finite map $$\pi_{n,k} :F_{n+k}\to {(F_{n+1})}^k_{F_n}.$$
To see this, look at the case
where $k=2$. If we consider an element $y\in F_{n+2}$, then $y$ can be
written as:
$$y=\{ (P_{1,1}, \ldots, P_{n,1}, P_{n+1,1}, P_{n+2,1}), \ldots,
(P_{1,d}, \ldots, P_{n,d}, P_{n+1,d}, P_{n+2,d})\}$$
The map $\pi_{n,2}$ is defined by sending $y$ to the following two elements of
$F_{n+1}$:
$$\{(P_{1,1}, \ldots, P_{n,1}, P_{n+1,1}),\ldots,
(P_{1,d}, \ldots, P_{n,d}, P_{n+1,d})\}$$ and
$$\{(P_{1,1}, \ldots, P_{n,1}, P_{n+2,1}), \ldots,
(P_{1,d}, \ldots, P_{n,d}, P_{n+2,d})\}$$ These two elements together
form an element of ${(F_{n+1})}^2_{F_n}$, thus defining $\pi_{n,2}(y)$.

Similiarly, an induction argument shows
that $F_{n+k}$ maps finitely to ${(F_{n+1})}^k_{F_n}$.

Notice that by definition the $E_n$'s and the $F_n$'s are not contained in
the {\it big diagonal} in $Y_n$; by the big diagonal, we mean the set
of $n$-tuples for which at least $2$ entries agree.  Consider the following
lemma from \cite{a}:

\begin{lemma}[see \cite{a}, Lemma 1]
Let $D\to Z$ be a generically finite morphism, and let
$\Delta_n$ be the big diagonal in $D_Z^n$.  Then there exists an integer
$n$ for which the $n$-th fiber product of $D$ over $Z$,
$D_Z^n\backslash\Delta_n\to Z$, is not dominant.
\end{lemma}

If we let $D=F_{n+1}$, $Z=F_n$, then this lemma shows that if $y\in
F_n$ then the dimension of the fiber above $y$ in $F_{n+1}$ is at
least $1$.  This means that any component of $F_n$ has a component of
$F_{n+1}$ above it of positive dimension. Moreover, the lemma shows
that if $y$ is in $F_n$ and if $k$ is a positive integer, then the
dimension of the fiber above $y$ in $F_{n+k}$ is at least $k$.

We want to perform an induction argument on the relative dimension of
$F_{n+1}$ over $F_n$ to show that $F_n$ is empty for all $n$. We know
that this dimension can't be greater that $d$, as this is the relative
dimension of $Y_{n+1}$ over $Y_n$, and we also know that it is at least $1$,
by the lemma above.  What we will do is show that, for any $l$, if the relative
dimension of $F_{n+1}$ over $F_n$ is at least $l$ then it must be at
least $l+1$. Eventually we will get to $l=d$ where the argument must
stop, and we will have to conclude that each $F_n$ is empty.

Assume that for all $n$ and for all $y\in F_n$, the dimension of the
fiber above $y$ in $F_{n+1}$ is greater than or equal to $l$.  Suppose
there exists an element $y\in F_n$ with the dimension of the fiber
above $y$ in $F_{n+1}$ exactly $l$.  By the semicontinuity of the
fiber dimension of projective maps, there is an irreducible component
of $F_n$, call it $M_n$, with the properties that $M_n$ contains $y$,
and the general fibers in $F_{n+1}$ over $M_n$ have dimension equal to
$l$. It also follows by induction that for any integer $k$ the
dimensions of the
fibers in $F_{n+k}$ above $M_n$ have dimension equal to $kl$.

Consider the following diagram where, as you may recall, the top map
is finite and birational.
$$
\begin{array}{ccc}
F_{n+k} & \longrightarrow & {(F_{n+1})}^k_{F_n}\\
        &                 &                   \\
        & \searrow        & \downarrow        \\
        &                 &                   \\
        &                 & F_n               \\
\end{array}
$$
Because the map $F_{n+k}\to {(F_{n+1})}^k_{F_n}$ is finite, there
exists at least one irreducible component $H_k$ of $F_{n+k}$
which dominates a component of ${(F_{n+1})}^k_{F_n}$; this
component of ${(F_{n+1})}^k_{F_n}$ dominates $M_n$ and has maximal relative
dimension $kl$ over $F_n$.

Later in this paper we will prove the following important proposition:

\begin{proposition}\label{prop} For large values of $k$ and $n$ every
component of the fiber product ${(F_{n+1})}^k_{F_n}$ of maximal
relative dimension is a variety of general type.
\end{proposition}

Therefore, for large $k$, $H_k$ dominates a variety of general type,
and because the map from $F_{n+k}$ to ${(F_{n+1})}^k_{F_n}$ is
birational, we can conclude that for large enough $k$ and $n$,
$F_{n+k}$ is also a variety of general type.  Remember, however, that
by definition the set of rational points in $F_{n+k}$ is dense.  Thus
Lang's conjecture combined with the work above implies a
contradiction.  Thus we have discomvered that for all $n$ and for all
$y\in F_n$, the dimension of the fiber above $y$ in $F_{n+1}$ must be
at least $l+1$.

Continuing by induction on this relative dimension, we will be forced
to conclude that if Lang's conjecture is true, we have a contradiction
unless $F_{n+k}$ is empty for large values of $k$ and $n$.  Since we
have surjective projections from $F_{n+1}\to F_n$ for all $n$,
however, this implies that all of the $F_n$'s must be empty, hence
proving Theorem \ref{mainth}.

The remainder of the paper will be devoted to proving Proposition
\ref{prop}.

\section{Background Lemmas and Definitions}

In definition \ref{def1}, we defined a big line bundle; now we provide an
extension of that definition which shall be quite useful.

\begin{definition} If $L$ is a line bundle on a variety $Z$ and $\cal J$ is
an ideal sheaf on $Z$, then we define $L\otimes \cal J$ to be { \em big} if
for high values of $k$, $H^0(L^{\otimes k}\otimes {\cal J}^k)$ induces
a birational map to projective space $P^N$ for some $N$.
\end{definition}

\begin{lemma}\label{curves.1} Assume $Z$ is a projective irreducible variety of
dim $l>0$, and that $$Z\subseteq C_1\times \cdots \times C_d$$ where the
$C_i$'s are irreducible projective curves.  Suppose further that for some
$i\in \{1,\ldots,d\}$ the projection map $\pi_i :Z\to C_i$ is surjective.
Then there exists a set $J\subseteq \{1, \ldots, d\}$ such that $i\in J$,
$\# J=l$, and the projection map $$\pi_J :Z\to \prod_{j\in J}C_j$$
is surjective.
\end{lemma}

{\bf Proof:} We proceed with induction on $d$.  If $d=1$, we simply
have $Z\subseteq C_1$; hence it must be that both $i=1$ and $l=1$.  We
use the fact that the inclusion map from $Z$ to $C_1$ must be either
constant or surjectve, but as $Z$ has positive dimension, it can't be
constant.  Therefore the lemma holds, with $J=\{1\}$.

Assume the result true for $d-1$ curves, and suppose that we have
$Z\subseteq C_1\times\cdots\times C_d$, with $Z$ surjecting onto $C_i$
for some $1\leq i\leq d$.  Choose $I\subseteq\{1, \ldots, d\}$
to be a subset of cardinality $d-1$ with the property that $i\in I$,
and let $\pi_I$ be the projection map $$\pi_I : Z\to\prod_{i\in I}C_i.$$
Define $Z_I=\pi_I (Z)$.  We know that $Z_I$ is
an irreducible projective variety, with dimension either $l$ or
$l-1$.

If $Z_I$ has dimension $l$, the induction hypothesis produces a
set $J\subseteq I$ such that $i\in J$,
$\# J=l$, and the projection map $$\pi_J :Z_I\to \prod_{j\in J}C_j$$
is surjective. Clearly $Z$
surjects to $Z_I$, so by composition we obtain a surjection
$$Z\to \prod_{j\in J}C_j.$$

Now suppose that $Z_I$ has dimension $l-1$.  The induction hypothesis
gives a set $J\subseteq I$ such that $i\in J$, $\# J=l-1$, and the
projection map $$\pi_J :Z_I\to \prod_{j\in J}C_j$$ is
surjective. Let $k\in \{1, \ldots, d\}$ be the one element which is not
in the set $I$.
Consider $Z_I\times C_k$. A simple dimension argument
shows that this is equal to $Z$.  Therefore we have a surjection
$$Z=Z_I\times C_k\to \prod_{j\in J}C_j\times C_k.$$
Therefore $Z$ surjects onto a product of $l$ curves, as desired.
Q.E.D.

\begin{lemma}\label{curves} Lemma \ref{curves.1} above still holds
if we replace
$Z\subseteq C_1\times\cdots\times C_d$ with a fiber product of
families of curves
$Z\subseteq C_1\times_B\cdots\times_B C_d$ where $B$ is a
projective variety, the $C_i$ form a family of curves over $B$, and that
for some $i$ the projection map from $Z$ the family $C_i$ is surjective.
\end{lemma}

{\bf Proof}:  Let $\eta\in B$ be a generic point.  Because $Z$ surjects to
each family $C_i$, it follows that each generic fiber of $Z$ surjects to
generic fibers of curves.  In other words, we may apply lemma
\ref{curves.1} to the fiber $Z_{\eta}$.  We obtain a set $J$ of size $l$
with $i\in J$ such that we have a surjective and generically finte
projection map
$$
\pi_J :Z_{\eta}\to \prod_{j\in J} C_{i,\eta}
$$
where the product above is a fiber product over $B$.  Because generic points
and fibers are dense, we know that the projection map
$$
\pi_J : Z\to \prod_{j\in J} C_i
$$
is generically finite (where the above product is once again a fiber
product over $B$).  The properness of the projection map implies the
surjectivity of $\pi_J$.  Q.E.D.

\begin{lemma}\label{surject} Suppose we have two surjective maps
$g_i:Z_i\to B$ for $i=1,2$.  Then the two natural maps
$f_i:Z_1 \times_B Z_2 \to Z_i$, $i=1,2$, are also surjective.
\end{lemma}

{\bf Proof:}  To show that $f_1$ is surjective, choose an element
$z_1\in Z_1$.  Let $b=g_1(z_1)$.  Then there exists a $z_2\in Z_2$
such that $g_2(z_2)=b$, as $g_2$ is surjective.  It follows that
$(z_1,z_2)\in Z_1\times_B Z_2$ and $f_1((z_1,z_2))=z_1$.  A similiar
argument shows that the map $f_2$ is also surjective. Q.E.D.

\begin{lemma}\label{product} Let $Y\subseteq Z_1\times Z_2$
be a variety, and let
$f_i: Y\to Z_i$ be surjective projection maps for $i=1,2.$
Suppose $L_1$ and $L_2$ are big line bundles on $Z_1$ and
$Z_2$ respectively.  Then the line bundle
$$
f_1^*L_1\otimes f_2^*L_2
$$
is big.  Further, if ${\cal J}_1$ and ${\cal J}_2$ are
ideal sheaves on $Z_1$ and $Z_2$, such that each $L_i\otimes {\cal J}_i$
is big on $Z_i$, then it follows that
$$
f_1^*(L_1)\otimes f_2^*(L_2)\otimes
(f_1^{-1}({\cal J}_1)\cdot f_2^{-1}({\cal J}_2))
$$
is big on $Y$.
Finally, the above still holds if we replace $Y\subseteq Z_1\times Z_2$
by a variety $Y$ which maps generically finitely to $Z_1\times Z_2$,
$$
h:Y\to h(Y)\subseteq Z_1\times Z_2
$$
provided that the projections from $h(Y)$ to $Z_i$ are still surjective.
\end{lemma}

{\bf Proof:} The first statement is a consequence of the second; hence
we just prove the second.

By assumption, since each $L_i\otimes {\cal J}_i$ is big on $Z_i$,
there exists open sets $U_i\subseteq Z_i$ such that for large $k$, global
sections of $H^0(Z_i, L_i^k\otimes {\cal J}^k)$ separate points in
$U_i$.  Let $U=f_1^{-1}(U_1)\cap f_2^{-1}(U_2)$, and let
$P=(P_1,P_2), Q=(Q_1,Q_2)\in U$.  We claim that we can produce
sections of
$$
H^0(Z_1\times Z_2, f_1^*(L_1)\otimes f_2^*(L_2)\otimes
(f_1^{-1}({\cal J}_1)\cdot f_2^{-1}({\cal J}_2)))
$$
for large $k$ which separate $P$ and $Q$.  In other words,
there is a section vanishing at $P$ and not at $Q$, and vice versa.
This will show that there are enough sections to induce a birational
map to projective space, hence this will suffice to show that
$$
f_1^*(L_1)\otimes f_2^*(L_2)\otimes
(f_1^{-1}({\cal J}_1)\cdot f_2^{-1}({\cal J}_2))
$$
is big.

We have the following map:
$$
f_1^*\otimes f_2^* :
H^0(L_1^{\otimes k}\otimes {{\cal J}_1}^k)\otimes
H^0(L_2^{\otimes k}\otimes {{\cal J}_2}^k)\to
H^0(f_1^*(L_1)\otimes f_2^*(L_2)\otimes
(f_1^{-1}({\cal J}_1)\cdot f_2^{-1}({\cal J}_2)))
$$
Therefore, sections of $H^0(Z_1\times Z_2,f_1^*(L_1)\otimes f_2^*(L_2)\otimes
(f_1^{-1}({\cal J}_1)\cdot f_2^{-1}({\cal J}_2)))$ have the form
$$
f_1^*(s_i)\otimes f_2^*(r_j)
$$
where $\{s_i\}$ and $\{r_j\}$ form bases for
$H^0(Z_1,L_1^{\otimes k}\otimes {\cal J}_1^k)$ and
$H^0(Z_2,L_2^{\otimes k}\otimes {\cal J}_2^k)$, respectively.

Because $P_1$ and $Q_1$ are elements of $U_1$, it follows that
there exists a section $s\in H^0(Z_1,L_1^{\otimes k}\otimes {\cal J}_1^k)$
such that $s(P_1)=0$ but $s(Q_1)\neq 0$.  Similiarly, there exists
a section $r\in H^0(Z_2,L_2^{\otimes k}\otimes {\cal J}^k_2)$
such that $r(P_2)=0$, but $r(Q_2)\neq 0$.  Consider the section
$f_1^*s\times f_2^*r$.  It follows that
$$
f_1^*s\times f_2^*r(P)=(s(P_1),r(P_2))=0
$$ and
$$
f_1^*s\times f_2^*r(Q)=(s(Q_1),r(Q_2))\neq 0.
$$
Similiarly, we can find a section of
$H^0(Z_1\times Z_2, f_1^*(L_1)\otimes f_2^*(L_2)\otimes
(f_1^{-1}({\cal J}_1)\cdot f_2^{-1}({\cal J}_2)))$,
which vanishes on $Q$ but not on $P$.  Therefore, for large $k$,
sections of
$H^0(Z_1\times Z_2, f_1^*(L_1)\otimes f_2^*(L_2)\otimes
(f_1^{-1}({\cal J}_1)\cdot f_2^{-1}({\cal J}_2)))$
generically separate points.  Hence, we
have that
$$
f_1^*(L_1)\otimes f_2^*(L_2)\otimes
(f_1^{-1}({\cal J}_1)\cdot f_2^{-1}({\cal J}_2))
$$
is big.

We now prove the last statement of the lemma.  Since
$h(Y)$ surjects both $Z_1$ and $Z_2$, the work above shows
that $f_1^*(L_1)\otimes f_2^*(L_2) \otimes (f_1^{-1}({\cal J}_1)\cdot
f_2^{-1}({\cal J}_2))$ is big on $h(Y)$.  Now
$h:Y\to h(Y)$ is a generically finite surjective morphism; hence the
pullback of a big sheaf
on $h(Y)$ will be big on $Y$. Q.E.D.

\section{The Proof}

Recall that for a given $n$, $F_n$ is contained in $Y_n={Sym}^d(X_B^n)$.
Let $X_n={(X_B^n)}^d$.  Then we have a map $\sigma : X_n\to Y_n$, namely the
quotient map given by action of $S_d$.  Now we can look at the inverse image
of $F_n$ under $\sigma$, which is contained in $X_n$.  Call it $G_n$.

These varieties $G_n$ will be central to our proof. We can think of the
$G_n$'s as forming a tower: we have projection maps
$\pi : G_{n+1} \to G_n$ such that
$$
\pi (G_{n+1})=G_n.
$$

Fix $n$.
Let $l$ be the relative dimension of $G_{n+1}$ over $G_n$.
Let $\phi$ be the projection of $G_n$ down to $B^d$, and for
$i=1, \ldots, d$, denote by $\pi_i$ the d projections from $B^d$ to $B$.
Notice that because the image of $G_n$ in ${(X_B^{n-1})}^d$ is $G_{n-1}$,
it follows that no matter which $n$ we start with the image
$\phi(G_n)$ in $B^d$ is the same.
Define
$B_i\subseteq B$ to be the image of an irreducible component of
$G_n$ in $B$ under the map
$\pi_i\circ \phi$. We shall see later that it is enough to restrict
our attention to one
component of $G_n$, hence, we do so now.

Choose $${\cal B}\subseteq B_1\times B_2\times
\cdots \times B_d$$ to be a component of $\phi(G_n)$
such that ${\cal B}$ surjects to each component
$B_i$.   Let $\widetilde{\cal B}$ be a desingularization of
$\cal B$.  By Hironaka's Theorem we know we can choose this desingularization
such that the discriminant locus is a divisor with normal crossings.

We have the following diagram:
$$\widetilde{\cal B}\to {\cal B}\subseteq B_1\times\ldots\times B_d
\subseteq B^d\to B$$
Notice that we have $d$ choices for the last map in the above diagram,
namely the $d$ projections $\pi_i$.  Let $g_i : \widetilde{\cal B}\to B$
be the compostion of the above where $\pi_i$ is chosen as the last mapping.

Define families of curves
$$
C_i=g_i^* (X\to B)=X\times_B \widetilde{\cal B}.
$$

$$
\begin{array}{ccc}
C_i 		   & \longrightarrow 		 & X \\
    		   &                 		 &   \\
\downarrow 	   &          			 & \downarrow \\
           	   &          			 &            \\
\widetilde{\cal B} & \stackrel{g_i}{\longrightarrow} & B \\
\end{array}
$$

We now pull back the varieties $G_n$ to $\widetilde{\cal B}$;
to simplify matters, we shall keep the same notation $G_n$.
We then also need to pull back the familes $C_i$ to $G_n$. Again, it
is easiest to abuse notation, and still refer to them as $C_i$.  In
our new situation we have the following containment:
$$G_{n+1}\subseteq C_1\times_{G_n} \cdots\times_{G_n} C_d.$$

Our goal is to work with these varieties $G_n$ to show that
for large $n$, $F_n$ is of general type, which will lead to our
contradiction. Our main difficulty will be in avoiding the singularities
of $G_n$.  To do this, we will eventually be working with products of
stable curves over the smooth base $\widetilde{{\cal B}}$.  These
curves will have canonical singularities, which are easier to control.
We would like to have $G_{n+1}$ surject onto each of these families
of curves $C_i$. While this is not always so, we will show that
we may reduce to this case.

For $i=1, \ldots d$, define projection maps $p_i :G_{n+1}\to C_i$.  To
illustrate the situation, assume for a moment that $G_{n+1}$ is
irreducible.  Then for each $i$, $p_i(G_{n+1})$ is an irreducible
subvariety of $C_i$.  Suppose that for one index $i$, the map is not
surjective; without loss of generality assume that $p_1(G_{n+1})=P$
where $P\subset C_i$ and has degree $k$ over $G_n$.  Now we go up $k$
steps and use the fact that
$$G_{n+k+1}\subseteq {(G_{n+1})}^{k+1}_{G_n}.$$ In other words,
$$G_{n+k+1}\subseteq P^{k+1}_{G_n}\times_{G_n}C_2^{k+1}\times_{G_n}
\cdots\times_{G_n} C_d^{k+1}.$$
Since $P$ consists of $k$ points over $G_n$, the pigeon-hole principle says
that two of the coordinates of $P^{k+1}_{G_n}$ must agree, contradicting the
fact that $G_{n+k}$ is not contained in the big diagonal.

Thus if $G_{n+1}$ is irreducible, we see that it is contained inside a product
of $d$ families of curves fibered over $G_n$ and that it surjects onto each
family.

Now suppose that $G_{n+1}$ consists of $m$ irreducible components:
$$
G_{n+1}=V_1\cup V_2\cup\cdots\cup V_m\subseteq C_1^n\times_{G_n} \ldots
\times_{G_n} C_d^n.
$$
Define projection maps $p_{i,j}:V_i\to C_j$ for $1\leq i\leq m$ and
$1\leq j \leq d$.  We want to show that at least one $V_i$ surjects onto all
$d$ curves.  So suppose for contradiction that
$$
p_{1,j(1)}(V_1)=S_1, \ldots, p_{m,j(m)}(V_m)=S_m
$$
where each $S_i\subset C_{j(i)}$ maps generically finitely to $G_n$ with
degree $r_i$.

Therefore, we know
$$
G_{n+1}\subseteq (p_{1,1}(V_1)\times_{G_n}\cdots\times_{G_n}p_{1,d}(V_1))
\cup\cdots\cup
(p_{m,1}(V_m)\times_{G_n}\cdots\times_{G_n}p_{m,d}(V_m))
$$
and
$$
G_{n+k}\subseteq {(G_{n+1})}^k_{G_n}
$$
So $G_{n+k}$ is contained in the above union fibered $k$ times over $G_n$.

In taking elements of
$$(p_{1,1}(V_1)\times_{G_n}\cdots\times_{G_n}p_{1,d}(V_1))$$
we have only $r_1$ distinct choices for the $j(1)$-st component.
Continuing in this spirit, in taking an element of
$$(p_{m,1}(V_m)\times_{G_n}\cdots\times_{G_n}p_{m,d}(V_m))$$
we have only $r_m$ distinct choices for the $j(m)$-th entry.

Thus, if we take $k\geq r_1+\cdots+r_m$, an element of $G_{n+k}$
must have $2$ coordinates equal, contradicting the fact that
$G_{n+k}$ is not contained in the big diagonal.

Therefore we know that if $G_{n+1}$ at least one
component of $G_{n+1}$ surjects onto all $d$ families of curves $C_i$.

Suppose once again that $G_{n+1}=V_1\cup\cdots\cup V_m$.  Define $W_1$ to
be the union of those $V_i$ which do not surject to each family $C_i$,
or are not dominant over $G_n$.  Define $W_2$ to be the union of those
components which do surject onto each $C_i$; by our work above, $W_2$ is
not empty.  Hence $$G_{n+1}\subseteq W_1\cup W_2$$

Choose a component $G$ of $G_{n+k}$ which is irreducible and dominant over
$G_n$ with relative dimension $kl$.  We then have
$$
G\subseteq G_{n+k}\subseteq {(G_{n+1})}^k_{G_n}\subseteq
{(W_1\cup W_2)}^k_{G_n}.
$$
If we let $J_1$ and $J_2$ run over all subsets of $\{1, \ldots, k\}$ such
that $J_1\cup J_2=\{1, \ldots, k\}$, then it follows that
$$
G\subseteq \cup_{J_1,J_2}{(W_1)}^{J_1}_{G_n}\times_{G_n}{(W_2)}^{J_2}_{G_n}
$$
Because $W_1$ and $W_2$ consist of unions of distinct components,
and $G$ is a component, it
follows that for one choice of $J_1$ and $J_2$, say $J_{1,k}$ and
$J_{2,k}$ that
$$
G\subseteq {(W_1)}^{J_{1,k}}_{G_n}\times_{G_n}{(W_2)}^{J_{2,k}}_{G_n}
$$
Recall that $G_{n+k}$ is not contained in the big diagonal. Therefore,
applying the same reasoning as earlier, we conclude that there exists
an integer $r$ such that $\#J_{1,k}\leq r$ for all $k$.  This is because,
by definition, the set $W_1$ consists of elements of $G_{n+1}$ which do
not surject to all of the families of curves $C_i$.  Since we know that
$\#J_{1,k}+\#J_{2,k}=k$, it follows that $\#J_{2,k}\geq k-r$, and so
as $k$ grows, the size of the sets $J_{2,k}$ also grows.  Since, by
definition, components of $W_2$ surject to each of the families $C_i$,
we can obtain a surjection
$$
{(W_2)}^{J_{2,k}}_{G_n}\to C_1^{k_1}\times_{G_n}\cdots\times_{G_n} C_d^{k_d}
$$
where we can make the exponents $k_i$ as large as we wish, by taking larger
values of $k$, since the size of $J_{2,k}$ grows with $k$.

Because $G$ is a component of $G_{n+1}$ of maximal dimension
which is dominant over $G_n$, it follows that
$G$ dominates a component $W$ of ${(W_2)}^{J_{2,k}}_{G_n}$.  This component
$W$ is contained inside a product of, say, $c$ of the components of $W_2$.
As $k$ grows, so does $k-r$; for large
values of $k-r$ at least one of those  $c$ components, call it $W'$,
will appear at
least $\kappa=\frac{k-r}{c}$ times, where $\kappa$ also grows with $k$.
Therefore $G$ dominates ${(W')}^\kappa_{G_n}$.

As $\kappa$ grows,  ${(W')}^\kappa_{G_n}$ will be mapping surjectively to
higher and higher powers of the families of curves $C_i$.  In the rest of
this paper we shall prove that
this is exactly what is needed to show that  ${(W')}^\kappa_{G_n}$ is a
variety of general type.  Moreover, we shall show that for $\kappa$ large
enough (i.e., for $k$ large enough) ${(W')}^{\kappa}_{G_n}$ modulo the
action of the symmetric group is a variety of general type.

Because $G$ dominates ${(W')}^\kappa_{G_n}$ we obtain a family of varieties
$$G\to {(W')}^\kappa_{G_n}.$$  Because $G_n$ is not contained in the fixed
point locus, over each general point in ${(W')}^\kappa_{G_n}$, the
fiber in $G$ consists of a product of curves of genus $g$ (as opposed to
a quotient of products of curves).  In other words,
each fiber is a variety of general type, as is the base.  We then utilize a
theorem of Viehweg (\cite{v}, Satz III)
to conclude that $G$ itself is a
variety of general type.  Viehweg's theorem states that if $Z\to B$ is
a family of varieties of general type where the base $B$ is also of
general type, then it follows that $Z$ is of general type.

Let $F$ be the image of $G$ in $F_{n+k}$.  Since $G$ is of maximal
dimension in $G_{n+k}$ it follows that $F$ is of maximal dimension
in $F_{n+k}$  Since ${(W')}^{\kappa}_{G_n}$ modulo the group
action is of general type, we may again apply Viehweg's theorem
as above to see that $F$ also is of general type, using that $F$
is a family of varieties of general type over the image of
${(W')}^{\kappa}_{G_n}$ modulo the group action.  Therefore, we've
shown that for large $n+k$, a component of $F_{n+k}$ of maximal
dimension is a variety of general type, proving proposition \ref{prop}.

Our work above proves the following proposition.

\begin{proposition}\label{prop.2}  In order to prove Proposition \ref{prop},
it suffices to prove it for the case where each $G_n$ is irreducible.
\end{proposition}

For the remainder of the paper,
we shall assume that for all large $n$ $G_{n+1}$ consists
of one irreducible component, which projects surjectively onto
each family of curves.

Because we are assuming that $G_{n+1}$ projects surjectively onto each
$C_i$, we may apply lemma \ref{curves}.  We thus obtain $d$
generically finite maps $\pi_i$ from $G_{n+1}$ to a product of the
families of curves $C_i$ fibered over $G_n$.  In other words, for each
$i=1,\ldots,d$ there exists a subset $J_i\subseteq \{1, \ldots, d\}$
such that $\#J_i=l$, $i\in J_i$, and
$$\pi_{J_i} :G_n\to \prod_{j\in J_i}C_j.$$

Recall, however, that we can project $G_n$ down to $\widetilde{\cal
B}$; this allows us to obtain $d$ generically finite surjective maps
from $G_{n+1}$ to a product of families of
curves fibered over $\widetilde{\cal B}$.

Since each $G_{n+k}$ is contained in a fiber
product of $G_{n+1}$ over $G_n$, we can apply the maps $\sigma_i$ to
each component of $G_{n+1}$ in this product, and we can produce a generically
finite map from $G_{n+k}$ to a product of powers of the $C_i$:
$$
G_{n+k}\to C_1^{k_1}\times_{\widetilde{\cal B}}\cdots
\times_{\widetilde{\cal B}} C_d^{k_d}.
$$
By taking $k$ large as necessary, we can increase the exponents
$k_i$, making them as large as we wish.

To ease notation, for large $n$ let
$$V_n=C_1^{k_1}\times_{\widetilde{\cal B}} \cdots
\times_{\widetilde{\cal B}} C_d^{k_d}.$$
So we have a generically finite map $G_n\to V_n$ where the exponents
appearing in $V_n$ can be made larger by increasing $n$.

Recall the statement of proposition \ref{prop}, which says that
for large values of $n$ and $k$, every component of
the fiber product ${(F_{n+1})}^k_{F_n}$ of maximal dimension is a variety
of general type.  We are ready to begin in earnest the proof of this
proposition.

First, we prove two more lemmas which will be of assistance to us.
The first statement in the following lemma is a well known fact, but
we include it here, as we will utilize it later.

\begin{lemma}\label{ideal} Let $f:V\to B$ be a flat morphism of irreducible
projective varieties with irreducible general fiber.  Let $L$ be a big line
bundle on $V$, and ${\cal J}\neq 0$ an ideal sheaf on $V$.
Let $$\pi_i:V_B^k\to V$$
be projections from the fiber product to the $i$-th factor.  Define
$L_k$ and ${\cal J}_k$ by
$$L_k=\otimes_i \pi_i^*L$$ and $${\cal J}_k=\sum_i\pi_i^{-1}{\cal J}.$$  Then:
\begin{enumerate}
\item There exists an integer $m$ such that $L^{\otimes m}\otimes {\cal J}$ is
big on $V$;
\item For
high enough values of $k$, it follows that $$L_k\otimes {\cal J}_k$$ is big
on $V_B^k$.
\end{enumerate}
\end{lemma}

{\bf Proof:} Consider the following exact sequence of cohomology groups:
$$
0\to H^0(V,L^{\otimes m}\otimes {\cal J})\to H^0(V,L^{\otimes m})\to
H^0(V,L^{\otimes m}/L^{\otimes m}\otimes {\cal J})\to 0
$$
Because $L$ is big,
we know that the dimension of $H^0(V,L^{\otimes m})$ is $O(m^{dim V})$.
Since the zero set of a variety has a lower dimension, we have that
the dimension of $H^0(V,L^{\otimes m}/L^{\otimes m}\otimes {\cal J})$ is less
than the dimension of $H^0(V,L^{\otimes m})$.  Therefore, the dimension
of $H^0(V,L^{\otimes m}/L^{\otimes m}\otimes {\cal J})$ is less than or equal
to $O(m^{dim V-1})$.  By counting the dimensions in the exact sequence,
we can conclude that the dimension of $H^0(V,L^{\otimes m}\otimes {\cal J})$ is
greater than or equal to $O(m^{dim V}-m^{dim V-1})$.  This growth in the
dimension of $H^0(V,L^{\otimes m}\otimes {\cal J})$ shows that there exists an
$m$ such that $L^{\otimes m}\otimes {\cal J}$ is big on $V$.

Because $\otimes \pi_i^*L^{\otimes m}=L_k^{\otimes m}$, we obtain a
map $$\otimes \pi_i^*H^0(V,L^{\otimes m})\to
H^0(V,\otimes\pi_i^*L^{\otimes m}).$$ The indicies $i$ range from
$1$ to $k$.  The inclusions
$$
\otimes \pi_i^* H^0(V,L^{\otimes m}\otimes {\cal J})
\hookrightarrow \otimes \pi_i^* H^0(V,L^{\otimes m})
$$
and
$$
H^0(V^k_B,L_k^{\otimes m}\otimes {\cal J}_k^k)\hookrightarrow
H^0(V^k_B,L_k^{\otimes m})
$$
produce a commutative diagram inducing a map
$$
\otimes \pi_i^* H^0(V,L^{\otimes m}\otimes {\cal J}) \hookrightarrow
H^0(V^k_B,{L_k}^{\otimes m}\otimes {\cal J}_k^k).
$$
Notice also that if we choose
$k>m$, we also have the inclusion
$$
H^0(V^k_B,L_k^{\otimes m}\otimes {\cal J}_k^k)\hookrightarrow
H^0(V^k_B,L_k^{\otimes m}\otimes {\cal J}_k^m)
$$
since ${\cal J}_k^k\subseteq {\cal J}_k^m$.

$$
\begin{array}{ccccc}
\otimes_i\pi_i^*H^0(L^{\otimes m}) & \longrightarrow
& H^0(L_k^{\otimes m}) &  &  \\
 & & & & \\
\uparrow &  & \uparrow & & \\
 & & & & \\
\otimes_i \pi_i^*H^0(L^{\otimes m}\otimes {\cal J}) & \longrightarrow &
H^0(L_k^{\otimes m}\otimes {\cal J}_k^k) & \longrightarrow &
H^0(L_k^{\otimes m}\otimes {\cal J}_k^m)\\
\end{array}
$$

As $L^{\otimes m}\otimes {\cal J}$ is big, we  have a birational map from
$V$ to $P^n$, induced by $H^0(V, L^{\otimes m}\otimes {\cal J})$.
The same reasoning as in Lemma \ref{product} lets us conclude that
we can use the sections of $H^0(V,L^{\otimes m}\otimes {\cal J})$ to
generically separate points of $V_B^k$.  Thus for large values of
$k$ it follows that $L_k\otimes {\cal J}_k$ is big on $V_B^k$.
Q.E.D.

\begin{lemma}\label{big} For large enough $n$, the relative dualizing sheaf
$$
\omega_{V_n/\widetilde{\cal B}}
$$
is big, and furthermore, the dualizing sheaf $$\omega_{V_n}$$ is also
big.
\end{lemma}

{\bf Proof:} Recall the following diagram:
$$\widetilde{\cal B}\to {\cal B} \hookrightarrow B_1\times \cdots \times
B_d\subseteq B^d\to B$$

Define $C_{i,0}$ to be the pullback of the family $X\to B$ to $B_i$.
Lemma 3.4 in \cite{chm} states that
the relative dualizing sheaf
$\omega_{C_{i,0}/B_i}$ is big.

Choose $n$ to be large, so that we have a generically finite map
$$
G_n\to C_1^{k_1}\times_{\widetilde{\cal B}}\cdots
\times_{\widetilde{\cal B}} C_d^{k_d}=V_n.
$$
Now for $i=1, \ldots, d$ consider the family
$$
{(C_{i,0})}_{B_i}^{k_i}\to B_i.
$$
Call
this map $f_i$.  Notice that
$$
{(C_{i,0})}_{B_i}^{k_i}\subseteq
{(C_{i,0})}^{k_i}.
$$
Now, $C_{i,0}$ maps surjectively to $B_i$ by
definition.  Hence lemma \ref{surject} implies that each projection
$$
{(C_{i,0})}_{B_i}^{k_i}\to C_{i,0}
$$
is surjective.  Thus lemma \ref{product}
tells us that since $\omega_{C_{i,0}/B_i}$ is big on $C_{i,0}$, it
follows that $\omega_{f_i}$ is big on ${(C_{i,0})}_{B_i}^{k_i}$.

Let $\phi$ be the map
$$
\phi : {\cal B} \hookrightarrow B_1\times \cdots \times B_d.
$$
Let $V_{n,0}$ be the pullback under $\phi$ of the family
$$
{(C_{1,0})}^{k_1}_{B_1}\times\cdots\times {(C_{d,0})}^{k_d}_{B_d}\to
B_1\times\cdots\times B_d
$$
Thus
$$
V_{n,0}\hookrightarrow {(C_{1,0})}^{k_1}_{B_1}\times\cdots\times
{(C_{d,0})}^{k_d}_{B_d}
$$
We also have surjective projections
$$
V_{n,0}\to {(C_{i,0})}^{k_i}_{B_i}
$$
They are surjective because
$\cal B$ was chosen so as to project onto each component of
$B_1\times\ldots\times B_d$.  Therefore, lemma \ref{product} tells us that
$\omega_{V_{n,0}/ {\cal B}}$
is big, since each $\omega_{f_i}$ is big on ${(C_{i,0})}^{k_i}_{B_i}$.

Let $\psi$ be the map
$$
\widetilde{\cal B}\to \cal B.
$$
We have that $V_n$ is the pullback of $V_{n,0}$ under $\psi$.  Thus
we know that
$$
\omega_{V_n/ \widetilde{\cal B}}=\psi^*\omega_{V_{n,0}/\cal B}
$$
and since $\psi$ is a generically finite map, it follows that
$\omega_{V_n/ \widetilde{\cal B}}$
is big.

It remains to show that the dualizing sheaf $\omega_{V_n}$ is big.
Recall that
$$
\omega_{V_n}=\omega_{V_n/ \widetilde{\cal B}}\otimes
\omega_{\widetilde{\cal B}}
$$
Choose $\cal I$ to be an ideal on the base $\widetilde{\cal B}$ such that
there is an injection $${\cal I}\hookrightarrow \omega_{\widetilde{\cal B}}.$$
It will be sufficient to show that
$$\omega_{V_n/\widetilde{\cal B}}\otimes {\cal I}
$$
is big.

To do this, define $Z$ to be the restriction of $X^d$ to
$\widetilde{\cal B}$.  In other words,
$$
Z=C_1\times_{\widetilde{\cal B}}\cdots\times_{\widetilde{\cal B}} C_d.
$$
By Lemma 3.4 in \cite{chm} we know that $\omega_{Z/\widetilde{\cal
B}}$ is big.  Recall that by choosing $n$ large enough we can control
the size of the exponents in $V_n$, making them as large as we wish.
Therefore, for any integer $m$ there exists an integer $n_m$ such that
if $n> n_m$, all the $k_i$ are greater than $m$.  Write
$$
k_i-m=r_i>0.
$$
Then we can rewrite
\begin{eqnarray*}
V_n & = & C_1^{r_1}\times_{\widetilde{\cal B}}\cdots
\times_{\widetilde{\cal B}}
C_d^{r_d}\times_{\widetilde{\cal B}} C_1^m\times_{\widetilde{\cal B}}
\cdots\times_{\widetilde{\cal B}}C_d^m\\
 & = & C_1^{r_1}\times_{\widetilde{\cal B}}\cdots
\times_{\widetilde{\cal B}}
C_d^{r_d}\times_{\widetilde{\cal B}} Z^m_{\widetilde{\cal B}}\\
 & = & V'\times_{\widetilde{\cal B}} Z^m_{\widetilde{\cal B}}\\
\end{eqnarray*}
where we define $V'$ to equal
$C_1^{r_1}\times_{\widetilde{\cal B}}\cdots
\times_{\widetilde{\cal B}}
C_d^{r_d}$.

The first part of this lemma shows that for large enough $n$,
$\omega_{V'/\widetilde{\cal B}}$ is big on $V'$.  Lemma
\ref{ideal} states that, since $\omega_{Z/ \widetilde{\cal B}}$ is big on
$Z$, for large enough $m$,
$$
\omega_{Z^m_{\widetilde{\cal B}}}\otimes {\cal I}$$ is big on
$Z^m_{\widetilde {\cal B}}$.

We now apply Lemma \ref{product} to see that
$$
\omega_{V'/\widetilde{\cal B}}\otimes
(\omega_{Z^m_{\widetilde{\cal B}}}\otimes {\cal I})$$
is big on $V'\times_{\widetilde{\cal B}} Z^m_{\widetilde{\cal B}}$.
In other words,
$$
\omega_{V_n/\widetilde{\cal B}}\otimes {\cal I}$$ is big, as desired;
hence $\omega_{V_n}$ is big on $V_n$.  Q.E.D.

We're almost ready to prove Proposition \ref{prop}.  Let us set up
some notation first.

Let $\widetilde{G_{n+k}}$ be an equivariant desingularization of $G_{n+k}$.
In other words, the action of a subgroup of the symmetric group $S_d$ on
$\widetilde{\cal B}$ and $G_{n+k}$ lifts to $\widetilde{G_{n+k}}$.
Let $r$ be the map from
$\widetilde{G_{n+k}}$ to $G_{n+k}$ which gives the resolution of
singularities.  Let $\Phi_{n+1}$ and $\Phi_{n+k}$ denote the set of points of
$G_{n+1}$ and $G_{n+k}$, respectively, which
are fixed by the group action.

Recall that we have a generically finite map, call it $\sigma_0$,
from $G_n$ to $V_n$:
$$
\sigma_0 :G_n\to V_n=C_1^{k_1}\times_{\widetilde{\cal B}} \cdots
\times_{\widetilde{\cal B}} C_d^{k_d}.
$$
For $i=1, \ldots, d$, denote by $\sigma_i$ the $d$ projections from $G_{n+1}$
to a product of $l$ families of curves, where $l$ is the relative dimension of
$G_{n+1}$ over $G_n$.  Let $Z_i$ be the image of $G_{n+1}$ under $\sigma_i$.
We can write
$$
Z_i=C_{j_1(i)}\times_{\widetilde{\cal B}}\cdots\times_{\widetilde{\cal B}}
C_{j_l(i)}
$$
where one of $j_1(i), \ldots, j_l(i)$ is equal to $i$.

Define varieties $V_{n+1,i}$ by
$$
V_{n+1,i}=V_n\times_{\widetilde{\cal B}}Z_i.
$$
Then for any positive integer $k$, we have
$$
V_{n+k,i}=V_n\times_{\widetilde{\cal B}}
{(Z_i)}^k_{\widetilde{\cal B}}.
$$
Let $\pi_0$ be the map from $G_{n+k}$ to $G_n$ given by projection to the
first $n$ coordinates, and
for $j=1, \ldots, k$, define projection maps
$$
\pi_j:G_{n+k}\to G_{n+1}
$$
as follows: If $(P_1, \ldots, P_n, \ldots, P_{n+k})\in G_{n+k}$ then
$$
\pi_j((P_1, \ldots, P_n, \ldots, P_{n+k})=(P_1, \ldots, P_n, P_{n+j}).
$$

Recall that $G_{n+k}\subseteq {(G_{n+1})}^k_{G_n}$.  We use this to define
maps
$$
(\sigma_0, \sigma_i, \ldots, \sigma_i):G_{n+k}\to V_{n+k,i}.
$$
They act in the following manner: $\sigma_0$ is applied to $G_n$, while the
$\sigma_i$ are applied to the copies of $G_{n+1}$.  In other words,
if $P=(P_1, \ldots, P_n, P_{n+1}, \ldots, P_{n+k})\in G_{n+k}$, then
$$(\sigma_0, \sigma_i, \ldots, \sigma_i)(P)=
(\sigma_0(P_1, \ldots, P_n), \sigma_i(P_{n+1}), \ldots, \sigma_i(P_{n+k})).
$$

We have the following diagram:

$$
\begin{array}{ccc}
 G_{n+k} & \stackrel{(\sigma_0, \sigma_i, \ldots, \sigma_i)}
{\longrightarrow} & V_{n+k,i}\\
 & & \\
\downarrow \pi_j & & \downarrow\\
 & & \\
G_{n+1} & \longrightarrow & V_{n+1,i}\\
 & & \\
\downarrow & & \downarrow\\
 & & \\
G_n & \stackrel{\sigma_0}{\longrightarrow} & V_n\\
\end{array}
$$

Recall from the proof of Lemma \ref{big} that we defined $Z$ to be the
restriction of $X^d$ to ${\widetilde{\cal B}}$.
$$
Z=C_1\times_{\widetilde{\cal B}}\cdots\times_{\widetilde{\cal B}}C_d.
$$
We know that $G_{n+1}$ maps to $Z$ and in fact, the $\sigma_i$'s
factor through this map. So let
$$
\tau :G_{n+1}\to W
$$
be the map from $G_{n+1}$ to $Z$, where $W$ is the image of $G_{n+1}$ in $Z$.
Let
$$
e:W\to Z
$$
be the inclusion map from $W$ to $Z$. Finally, define the projection maps
$$
\rho_i :Z\to Z_i
$$
where
$$
\rho_i (C_1\times_{\widetilde{\cal B}}\ldots\times_{\widetilde{\cal B}}C_d)
=C_{j_1(i)}\times_{\widetilde{\cal B}}\ldots\times_{\widetilde{\cal B}}
C_{j_l(i)}.
$$

We have the diagram:

$$
\begin{array}{ccccl}
 & & W & & \\
 & \nearrow & & \searrow e & \\
 & & & & \\
G_{n+1} & & \longrightarrow & & Z\\
 & & & & \\
 & & \searrow & & \downarrow \rho_i\\
 & & & & \\
 & & & & Z_i\\
\end{array}
$$

An important fact to notice is that the image in $W$ of the fixed points
$\Phi_{n+1}$ is not all of $W$.  Why is this?  Well, let's look at what the
map $\tau :G_{n+1}\to W\subseteq Z$ does.  Remember that
$G_{n+1}\subseteq {(X_B^{n+1})}^d$ and that $Z=X^d|_{\widetilde{\cal B}}$.
We can write a $K$-rational element, $g$, of $G_{n+1}$ as
$$
g=({(P_1, \ldots, P_{n+1})}^{\alpha_1}, \ldots,
{(P_1, \ldots, P_{n+1})}^{\alpha_d})
$$
where $\alpha_1, \ldots, \alpha_d$ are the embeddings of $L$ into
$\bar{K}$ fixing $K$.  Then
$$
\tau (g)=({P_{n+1}}^{\alpha_1}, \ldots, {P_{n+1}}^{\alpha_d}).
$$
Therefore $\tau (g)$ is actually a Galois orbit.  Recall that the
$P_{n+1}$'s are defined over $L$ but not any smaller field between $L$
and $K$.  Hence all of ${P_{n+1}}^{\alpha_1}, \ldots, {P_{n+1}}^{\alpha_d}$
are distinct.  Therefore, $\tau (G_{n+1})=W$ is the closure of a set of points
which are not fixed by the group action.

This fact that the fixed points do not surject to $W$ will soon be very
important.

We want to prove that every component of $F_{n+k}$ of maximal dimension is
a variety of general type.  This means that if $\widetilde{F_{n+k}}$ is
a resolution of singularities of $F_{n+k}$, we need to show that
the dualizing sheaf $\omega_{\widetilde{F_{n+k}}}$ is big.

The main idea is to pull back sections of dualizing sheaves along our
various projections. We know that both $\omega_{V_n}$ and
$\omega_{V_{n+k,i}}$ are big for large $n$.  ($\omega_{V_{n+k,i}}$
is big since the families of the curves $C_i$ appear in $V_{n+k,i}$ with
even larger exponents.)  We would like to utilize the fact that
$\omega_{V_{n+k,i}}$ is big to pull back a lot of sections to
$\widetilde{G_{n+k}}$ which vanish along the fixed points of the group
action (for we will see that this is what suffices for $F_{n+k}$ to
be of general type).  Unfortunately we can't pull back sections along
any one projection of $\widetilde{G_{n+k}}$ to $V_{n+k,i}$; the
sections pulled back in this way will not vanish along the fixed
points.  To overcome this difficulty we will pull back sections via
all the various projections, tensor them together, and use our earlier
lemmas to obtain that $\omega_{\widetilde{G_{n+k}}}$ is big, with
lots of sections vanishing on the fixed points.

Recall that $G_{n+k}$ is contained in a fiber power of families of
stable curves over $\widetilde{\cal B}$, and, moreover,
$\widetilde{\cal B}$ was chosen to be smooth and irreducible, with the
property that its discriminant locus is a divisor of normal crossings.
We appeal to Lemma 3.3 of \cite{chm}, which states that this implies
that any singularities of $V_{n+k,i}$ are canonical.  This means that
the singularities do not prohibit the extension of pluricanonical
sections from the smooth locus to any desingularization.  Therefore we
may pull back sections of $V_{n+k,i}$ to $\widetilde{G_{n+k}}$.  We
have the following injection:
$$
r^*(\sigma_0,\sigma_i, \ldots,\sigma_i)^*
\omega_{V_{n+k,i}}\hookrightarrow
\omega_{\widetilde{G_{n+k}}}
$$
and so
$$
r^*(\otimes_i((\sigma_0,\sigma_i, \ldots,\sigma_i))^*
\omega_{V_{n+k,i}}\hookrightarrow
(\omega_{\widetilde{G_{n+k}}})^{\otimes d}
$$
Now, the left-hand side above can be rewritten as:
$$
r^*(\pi_0^*\sigma_0^*\omega_{V_n})^{\otimes d}\otimes
r^*(\pi_1^*\otimes_i\sigma_i^*\omega_{Z_i/\widetilde{\cal B}})
\otimes\ldots\otimes
r^*(\pi_k^*\otimes_i\sigma_i^*\omega_{Z_i/\widetilde{\cal B}})
$$
Notice in the second diagram above that
$$
\sigma_i=\rho_i\circ e\circ \tau
$$
Hence,
$$
\sigma_i^*=\tau^* e^* \rho_i^*.
$$
Let's look at one of the last $k$ terms above (take some $j$ such
that $1\leq j\leq k$):
\begin{eqnarray*}
r^*(\pi_j^*\otimes_i\sigma_i^*\omega_{Z_i/\widetilde{\cal B}}) & = &
r^*(\pi_j^*(\otimes_i\tau^* e^* \rho_i^*\omega_{Z_i/\widetilde{\cal B}}))\\
 & = & r^*(\pi_j^*(\otimes_i\tau^* e^*
(\omega_{C_{j_1(i)}/\widetilde{\cal B}}\otimes\cdots\otimes
\omega_{C_{j_l(i)}/\widetilde{\cal B}})))\\
 & = & r^*(\pi_j^*\tau^* e^*(\otimes_i
(\omega_{C_{j_1(i)}/\widetilde{\cal B}}\otimes\cdots\otimes
\omega_{C_{j_l(i)}/\widetilde{\cal B}})))\\
 & = & r^*(\pi_j^*\tau^*e^*(\omega_{C_1/\widetilde{\cal B}}^{\otimes l_1}
\otimes\cdots\otimes\omega_{C_d/\widetilde{\cal B}}^{\otimes l_d}))\\
\end{eqnarray*}
The exponents $l_i$ are defined by the last line above, and we are
guaranteed that each one is positive.  We know that for each $i$,
one of the $j_k(i)$'s is equal to $i$; hence we when we tensor over all $i$
we know get each family $C_i$ appearing some positive number $l_i$ times.

Let $M$ denote the line bundle
$$
M=\omega_{C_1/\widetilde{\cal B}}^{\otimes l_1}
\otimes\cdots\otimes\omega_{C_d/\widetilde{\cal B}}^{\otimes l_d}.
$$
We claim that $M$ is a big line bundle on $Z$, and if we let $M_W$
denote $e^*M$, then, in fact, $M_W$ is big on $W$.

Recall that
$$
M=\otimes_i \rho_i^*\omega_{Z_i/\widetilde{\cal B}}.
$$
In other words,
$$
M=\omega_{C_1/\widetilde{\cal B}}^{\otimes l_1}\otimes\cdots\otimes
\omega_{C_d/\widetilde{\cal B}}^{\otimes l_d}
$$
where each $l_i$ is positive.

As in the proof of Lemma \ref{big}, let $C_{i,o}$ be the pullback of
$X\to B$ to $B_i\subseteq B$; recall that $\omega_{C_{i,0}/B_i}$ is big.
Moreover, using techniques from the proof of Lemma \ref{big}, we know
that $\omega_{{(C_{i,0})}^{l_i}_{B_i}}$ is big.

We have
$$
Z\subseteq C_{1,0}\times\cdots\times C_{d,0}\to B_1\times\cdots\times B_d
$$
and each $Z$ surjects to each $B_i$.  Since each
$\omega_{C_{i,0}/B_i}$ is big, it follows that each
$\omega^{\otimes l_i}_{C_{i,0}/B_i}$ is also big.  Therefore
applying Lemma \ref{product} we see that the pullback
of $\otimes_i\omega^{\otimes l_i}_{C_{i,0}/B_i}$ to $Z$ is big.  This
pullback is equal to
$$\omega_{C_1/\widetilde{\cal B}}^{\otimes l_1}\otimes\cdots\otimes
\omega_{C_d/\widetilde{\cal B}}^{\otimes l_d}=M$$
so $M$ is big on $Z$.  Now $W\subseteq Z$ and $W$ surjects to each
$C_i$ since $W$ is the image of $G_{n+1}$. Hence applying lemma
\ref{product} once again, we conclude that $M_W$ is big on $W$.

Putting all of the above together, we see that
$$
r^*(\pi_0^*\sigma_0^*\omega_{V_n})^d\otimes r^*(\pi_1^*\tau^* M_W)
\otimes\cdots\otimes r^*(\pi_k^*\tau_* M_W)\hookrightarrow
{\omega_{\widetilde{G_{n+k}}}}^{\otimes d}.
$$
Recall that the image of the fixed points $\Phi_{n+1}$ in $W$ is not all
of $W$.  Therefore this image forms a proper subvariey
of $W$, and we obtain a non-zero ideal $\cal I$ of sections of
$M_W$ which vanish on this subvariety.

We now apply the following lemma, borrowed from
\cite{chm}.

\begin{lemma}\label{fixed}[see \cite{chm}, Lemma 4.1]
Suppose that $X$ is a projective variety with
$G$ a finite group of order $l$ acting on $X$.  Let $\theta$ be an $m$
canonical form on $X$ which is invariant under the group action of $G$.
Then if $\theta$ vanishes to order $m(l-1)$ on the locus of all points
of $X$ fixed by the group action, then $\theta$ descends to smooth form
on $X/G$.
\end{lemma}

We use Lemma \ref{fixed} to show that
smooth pluricanonical forms on $\widetilde{G_{n+k}}$ descend to smooth
forms on $\widetilde{F_{n+k}}$ modulo the group action, if the forms vanish
to a prescribed order on the fixed point locus.  Bear in mind that
$\widetilde{G_{n+k}}$ modulo the group action is nothing more than
$\widetilde{F_{n+k}}$, a desingularization of $F_{n+k}$.  So to show that
$\omega_{\widetilde{F_{n+k}}}$ is big for large $k$, we will show that
$\omega_{\widetilde{G_{n+k}}}$ not only is big, but also has lots of
sections vanishing to arbitrarily high order on the fixed point locus.

By the first statement of Lemma \ref{ideal}, we know that for some $s$
the line bundle $M_W^{\otimes s}\otimes
{\cal I}$ is big.

Look at the following diagram:
$$
\begin{array}{ccccc}
\widetilde{G_{n+k}} & {r\atop\longrightarrow} & G_{n+k} &
\hookrightarrow & G_n\times W\times\cdots\times W \\
 & & & & \\
 & & & \searrow & \downarrow \\
 & & & & \\
 & & & & V_n\times W\times\cdots\times W \\
\end{array}
$$
Let ${\cal J}_\Phi$ be the locus of fixed points in $G_{n+k}$.
Then
$$
\sum_j \pi_j^{-1}\tau_j^{-1}{\cal J}\subset {\cal J}_\Phi
$$
and so
$$
\prod_j \pi_j^{-1}\tau_j^{-1}{\cal J}\subset{\cal J}^k_{\Phi}
$$
Therefore,
$$
r^*{(\pi_0^*\sigma_0^*\omega_{V_n})}^{\otimes ds}
\otimes r^*(\otimes_j\pi_j^*\tau^*M_W^{\otimes s}
\otimes(\prod_j  \pi_j^{-1}\tau_j^{-1}{\cal J}))\hookrightarrow
\omega_{\widetilde{G_{n+k}}}^{\otimes ds} \otimes {\cal J}^k_\Phi
$$

We know that $\omega_{V_n}$ is big, and
as in the proof of lemma \ref{ideal}, as long as $k>s$,
$M_W^{\otimes s}
\otimes(\prod_j  \pi_j^{-1}\tau_j^{-1}{\cal J}))$ is big.
Therefore, for $k>s$ the entire left hand side above is
big, so we have lots of sections of
$\omega_{\widetilde{G_{n+k}}}^{\otimes ds}$
vanishing to high order, and the proof is complete.

 \normalsize \vspace*{1 cm}


\begin{thebibliography}{HHHHHHH}

\bibitem[A]{a} D. Abramovich, {\em Uniformit\'e des points
rationnels des courbes alg\'ebriques sur les extensions quadratiques
et cubiques}, C.R. Acad. Sci. Paris, t. 321, S\'erie I,
p. 755-758, 1995.\\

\bibitem[CHM]{chm} L. Caporaso, J. Harris, B. Mazur, {\em Uniformity
of rational points}, J. Amer. Math. Soc., to appear\\
{\tt ftp://ftp.math.harvard.edu/pub/uniformityofrationalpoints.tex}

\bibitem[V]{v} E. Viehweg, {\em Die additivit\"at der Kodaira
dimension f\"ur projektive Fasserr\"aume \"uber variet\"aten des
allgemeinen typs,}, Jour. reine und angew. Math. 330, p. 132-142, 1982.

\end{thebibliography}
\end{document}